\documentclass[preprint,12pt,3p]{elsarticle}

\usepackage{amssymb}
\usepackage{amsmath,amsthm}
\usepackage{mathrsfs}
\usepackage{hyperref}
\usepackage{cleveref}
\usepackage[all,cmtip]{xy}
\usepackage{epsf, graphicx}
\usepackage{latexsym,amsfonts,amsbsy,amssymb}
\usepackage{geometry}
\usepackage{titletoc}
\usepackage{rotating}
\usepackage{algorithm}
\usepackage{algorithmic}
\usepackage{amscd}

\makeatletter

\@addtoreset{equation}{section} \makeatother
\newtheorem{theorem}{Theorem}[section]
\newtheorem{lemma}[theorem]{Lemma}
\newtheorem{corollary}[theorem]{Corollary}

\def\O{\mathcal{O}}

\def\ba{{\rm bar}}

\makeatletter
\def\ps@pprintTitle{%
\let\@oddhead\@empty
\let\@evenhead\@empty
\def\@oddfoot{\reset@font\hfil\thepage\hfil}
\let\@evenfoot\@oddfoot
}
\makeatother

\begin{document}

\begin{frontmatter}

\title{Descent of splendid Rickard equivalences in alternating groups}

\author{Xin Huang}
\address{SICM, Southern University of Science and Technology, Shenzhen 518055, China}

\begin{abstract}

We show that each block of an alternating group over an arbitrary complete discrete valuation ring is splendidly Rickard equivalent to its Brauer correspondent. This provides new evidence for a refined version of Brou\'{e}'s abelian defect group conjecture proposed by Kessar and Linckelmann.

\end{abstract}


\end{frontmatter}


\section{Introduction}\label{s1}

Brou\'{e}'s abelian defect group conjecture is an important conjecture in the representation theory of finite groups. In \cite{Kessar_Linckelmann}, Kessar and Linckelmann proposed a refined version of the abelian defect group conjecture, namely that for any complete discrete valuation ring $\O$ and any block of a finite group over $\O$ with an abelian defect group, there is a splendid Rickard equivalence between the block algebra and its Brauer correspondent. In Brou\'{e}'s original conjecture, the complete discrete valuation rings are assumed to be with splitting residue fields. Discovering more evidence for the refined conjecture is one of many motivations for realising known Rickard equivalences over non-splitting fields.

Recently, Boltje (\cite[Theorem 1.4]{Boltje}) proved that, roughly speaking, if a splendid Rickard equivalence between a block and its Brauer correspondent over a sufficiently large discrete valuation ring descends to the ring of $p$-adic integers, then various refinements of the Alperin-McKay conjecture, proposed by Isaacs-Navarro (\cite{IN}), Navarro (\cite{Navarro}) and Turull (\cite{T08},\cite{T13}) hold for this block.
This is another important reason for pursuing the question of descent of splendid Rickard equivalences.

Brou\'{e}'s original conjecture for alternating groups has been proved by Marcus \cite{Marcus}. In this note, we prove that the refined conjecture holds for alternating groups.

Throughout this note, $p$ is a prime number, $k'$ is a perfect field of characteristic $p$, $\O'$ is a complete discrete valuation ring of characteristic $0$ with residue field $k'$. Let $k$ be the prime field of $k'$, there is a unique complete discrete valuation ring $\O$ contained in $\O'$ such that $J(\O)=p\O$ and that the image of $\O$ under the canonical surjection $\O'\to k'$ is $k$ (see \cite[Chapter 2, Theorem 3, 4 and Proposition 1]{Serre}, or the paragraph proceeding \cite[Definition 1.8]{Kessar_Linckelmann}). Clearly $\O\cong \mathbb{Z}_p$, the ring of $p$-adic integers.

Let $G$ be a finite group. By a {\it block} of the group algebra $\O'G$, we mean a primitive idempotent $b$ of the center of $\O'G$, and $\O'Gb$ is called a {\it block algebra}. Let $H$ be another finite group, and let $c$ be a block $\O'H$. The block algebras $\O'Gb$ and $\O'N_G(P)c$ are called {\it Rickard equivalent} if there is a bounded complex $X$ of $(\O'Gb,\O'Hc)$-bimodules, each projective as a left $\O'Gb$-module and as a right $\O'Hc$-module, such that $X\otimes_{\O'Hc} X^*$ and $\O'Gb$ are homotopy equivalent as complexes of $(\O'Gb,\O'Gb)$-bimodules, and $X^*\otimes_{\O'Gb} X$ and $\O'Hc$ are homotopy equivalent as complexes of $(\O'Hc,\O'Hc)$-bimodules, where $X^*$ denotes the $\O'$-dual of $X$. The complex $X$ is called a {\it Rickard complex}. Suppose that $\O'Gb$ and $\O'Hc$ have a common defect group $P$, and each term of $X$, considered as an $\O'(G\times H)$-module, is a $p$-permutation module whose indecomposable summands have vertices contained in $\{(u,u)|u\in P\}$. Then $X$ is called {\it splendid}, and $\O'Gb$ and $\O'Hc$ are called {\it splendidly Rickard equivalent}.


Let $n>1$ be a positive integer. It is well-known that $k$ is already a splitting field of the symmetric group $S_n$. Chuang and Rouquier (\cite[Theorem 7.6]{CR}) proved that any block algebra of $\O S_n$ with a non-trivial abelian defect group is splendidly Rickard equivalent to its Brauer correspondent.

Under the assumption that $k'$ is algebraically closed, in \cite{Marcus}, Marcus proved the existence of a splendid Rickard equivalence between any block algebra of $\O' A_n$ with abelian defect group and their Brauer correspondent algebra by using the Chuang-Rouquier Theorem and a descent result (\cite[Theorem 4.7]{Marcuscomm}). Using some known results, it is easy to show that when $p>2$, any block of $\O'A_n$ with a non-trivial defect group is contained in $\O A_n$ (see Lemma \ref{block} below). Based on Marcus' work, we prove the following theorem, which implies the refined abelian defect group conjecture holds for alternating groups.


\begin{theorem}\label{main}
Each block algebra of $\O A_n$ with a non-trivial abelian defect group is splendidly Rickard equivalent to its Brauer correspondent algebra.
\end{theorem}

Our way to prove Theorem \ref{main} (in Section \ref{proof}) is checking that Marcus' proof can be modified into a proof for the ring $\O$, using some descent results in \cite{Kessar_Linckelmann}. Theorem \ref{main} together with \cite[Theorem 1.4 (c)]{Boltje} imply the following corollary.

\begin{corollary}
Turull's conjecture (\cite[Conjecture]{T13}) holds for blocks of alternating groups with abelian defect groups.
\end{corollary}




\section{Preliminaries}\label{notation}

In this section, we briefly review some material in \cite[Section 2]{Marcus}. Unless specified otherwise, modules in this note are left modules.

\subsection{Graded algebras and modules} Let $C_2:=\langle \sigma\rangle$ be a cyclic group of order $2$. Let $\Lambda\in\{\O',k'\}$. Assume that $p>2$.  The group $\hat{C}_2:={\rm Hom}(C_2,{\Lambda}^\times)$ is isomorphic to $C_2$, and we have $\hat{C}_2=\langle\hat{\sigma}\rangle$, where $\hat{\sigma}(\sigma)=-1$.

Let $R=R_1\oplus R_\sigma$ be a $C_2$-graded $\Lambda$-algebra, not necessary strongly graded.
Then $\hat{C}_2$ acts on $R$ via $C_2$-graded algebra automorphisms by
\begin{equation}\label{action}
{}^{\hat{\rho}}r_g=\hat{\rho}(g)r_g,
\end{equation}
for all $g\in C_2$ and $\hat{\rho}\in \hat{C}_2$. Then components of $R$ can be recovered as
$$R_{\sigma^j}=\{r\in R~|~{}^{\hat{\sigma}}r=(-1)^jr\},$$
for $j=0,1$. Define the {\it skew group algebra} $R\rtimes \hat{C}_2$: let $R\rtimes \hat{C}_2$ be equal to $R\otimes_\Lambda \Lambda\hat{C}_2$ as an $\Lambda$-module, with multiplication defined by
$$(r\otimes \hat{\rho})(r'\otimes \hat{\rho}')=r\cdot {}^{\hat{\rho}}r'\otimes \hat{\rho}\hat{\rho}',$$
where $r,r'\in R$ and $\hat{\rho},\hat{\rho}'\in\hat{C}_2$. Similarly, we can define $\hat{C}_2\ltimes R$.

Let $M=M_1\oplus M_\sigma$ be a $C_2$-graded $R$-module. Then $M$ becomes an $R\rtimes \hat{C}_2$-module with the action defined by
$$(r\otimes \hat{\rho})m_g=\hat{\rho}(g)rm_g$$
for all $r\in R$, $g\in C_2$, $m_g\in M_g$ and $\hat{\rho}\in \hat{C}_2$. If $f:M\to M'$ is a homomorphism of $C_2$-graded $R$-modules, then one easily checks that $f$ is also a homomorphism of $R\rtimes \hat{C}_2$-module. Conversely, let $M$ be an $R\rtimes \hat{C}_2$-module, and let
$$M_{\sigma^j}:=\{m\in M~|~\hat{\sigma}m=(-1)^jm\}.$$
Then every element $m\in M$ can be written uniquely as $m=m_1+m_{\sigma}$, with $m_{\sigma^j}\in M_{\sigma^j}$ for $j=0,1$, where $m_1=(m+\hat{\sigma} m)/2$, $m_{\sigma}=(m-\hat{\sigma}m)/2$. It is a routine exercise to check that $R_gM_h\subseteq M_{gh}$ for any $g,h\in C_2$, and if $f:M\to M'$ is a homomorphism of $R\rtimes \hat{C}_2$-modules, then $f(M_g)\subseteq M'_g$ for any $g\in C_2$. So, the following lemma holds.

\begin{lemma}[{\cite[Proposition 2.3]{Marcus}}]\label{gradcategory}
The category of $C_2$-graded $R$-modules is isomorphic to to category of $R\rtimes \hat{C}_2$-modules. 
\end{lemma}

Let $R$ and $S$ be two $C_2$-graded $\Lambda$-algebras. Denote by $S^{\rm op}$ the opposite algebra of $S$. Then $\hat{C}_2$ acts on $R\otimes_{\Lambda}S^{\rm op}$ diagonally by $${}^{\hat{\rho}}(r\otimes s)={}^{\hat{\rho}}r\otimes {}^{\hat{\rho}^{-1}}s.$$ So we can consider the skew group algebra $(R\otimes_{\Lambda}S^{\rm op})\rtimes \hat{C}_2$. As Lemma \ref{gradcategory}, the category of $C_2$-graded $(R,S)$-bimodules is isomorphic to the category of $(R\otimes_{\Lambda}S^{\rm op})\rtimes \hat{C}_2$-modules.

Let $M$ be an $(R,S)$-bimodule and $\rho\in \hat{C}_2$, define an $(R,S)$-bimodule ${}^{\hat{\rho}}M$:  let ${}^{\hat{\rho}}M=M$ as an $\Lambda$-module with the action
\begin{equation}\label{2.1}
(r\otimes s)\cdot_{\hat{\rho}} m:=({}^{\hat{\rho}^{-1}}r\otimes {}^{\hat{\rho}}s)m
\end{equation}
for all $m\in M$, $r\in R$, $s\in S$ and $\hat{\rho}\in\hat{C}_2$.

\subsection{Graded algebras and Galois automorphisms}

Let $R$ be an $k$-algebra, and set $R':=k'\otimes_kR$, an $k'$-algebra. Let $\Gamma:={\rm Gal}(k'/k)$. For an $R'$-module
$U$ and a $\tau\in \Gamma$, denote by ${}^\tau U$ the $R'$-module which is equal to $U$ as a module over the subalgebra $1\otimes R$ of $R'$, such that $\lambda\otimes r$ acts on $U$ as $\tau^{-1}(\lambda) \otimes r$ for all $r \in R$ and $\lambda\in k'$.
The $R'$-module $U$ is {\it $\Gamma$-stable} if ${}^\tau U\cong U$ for all $\tau\in \Gamma$.

Let $G$ be a finite group. Assume further that $R'$ is a $G$-graded $k'$-algebra via the decomposition $R'=\oplus_{g\in G}R'_g$, and that $U=\oplus_{g\in G}U_g$ is a $G$-graded $R'$-module. Then ${}^\tau U$ is also a $G$-graded $R'$-module via the same decomposition ${}^\tau U=\oplus_{g\in G}U_g$ of ${}^\tau U$ as a $k'$-module.

\subsection{Blocks of $\O'S_n$ and $\O'A_n$}

Let $n>1$ be a positive integer. Let $G:=S_n$, and $G^+:=A_n$. Assume that $p>2$, and that $k'$ is a finite splitting field for all subgroups of $G$. As an $\O'$-module, the group algebra $R:=\O'G$ can be decomposed as $R_1\oplus R_\sigma$, where $R_1=\O' G^+$, $R_\sigma=\O'G^+t$, and $t=(1,2)$, an element of $G\backslash G^+$. With respect to this decomposition, $R$ is a $C_2$-graded $\O'$-algebra.
Let $b^+$ be a block of $\O'G^+$ with a nontrivial defect group $P$. Let $H^+:=N_{G^+}(P)$, and $c^+\in \O'H^+$ be the Brauer correspondent of $b^+$. Let $b$ be a block of $\O'G$ covering $b^+$, then $P$ is also a defect group of $b$. Let $H:=N_G(P)$, $c\in \O'H$ be the Brauer correspondent of $b$. Then by the Harris-Kn\"{o}rr correspondence, $c$ covers $c^+$.

By an analog of Frattini's argument, we have $HG^+=G$ and hence $G/G^+\cong H/H^+\cong C_2$. Hence $S:=\O'H$ is a $C_2$-graded $\O'$-algebra. $\hat{C}_2$ acts on $R$ and $S$ via automorphisms of $C_2$-graded algebras. The Brauer correspondent of ${}^{\hat{\sigma}}b$ is ${}^{\hat{\sigma}}c$. Denote by $\hat{C}_{2,b}$ the stabiliser of $b$ under this action. The conjugation of $G$ (resp. $H$) on $G^+$ (resp. $H^+$) induces an action of $C_2$ on blocks of $\O' G^+$ (resp. $\O'H^+$). Denote by $C_{2,b^+}$ the stabiliser of $b^+$ in $C_2$. Then by \cite[Lemma 9.9]{Dade}, we have that
$$\sum_{g\in[C_2/C_{2,b^+}]}{}^gb^+
=\sum_{\hat{\rho}\in[\hat{C}_2/\hat{C}_{2,b}]}{}^{\hat{\rho}}b$$
as elements of $\O'G$ (also $\O'G^+$), where $[C_2/C_{2,b^+}]$ denotes a full set of representatives for the left cosets of $C_{2,b^+}$ in $C_2$. By the first paragraph of \cite[\S 5]{FH} and the first paragraph of \cite[\S 7]{FH}, we can deduce that $C_{2,b^+}=C_2$. So
\begin{equation}\label{2.1}
b={}^{\hat{\sigma}}b=b^+~~~{\rm or}~~~b\neq{}^{\hat{\sigma}}b,~b^+=b+{}^{\hat{\sigma}}b.
\end{equation}
Since both $b$ and ${}^{\hat{\sigma}}b$ are contained in $\O G$, we see that $b^+\in \O G^+$. We write this as a lemma.

\begin{lemma}\label{block}
Let $P$ be a nontrivial $p$-subgroup of $G^+$. An element of $\O'G^+$ is a block of $\O'G^+$ having $P$ as a defect group if and only if it is a block of $\O G^+$ having $P$ as a defect group.
\end{lemma}


\section{Proof of Theorem \ref{main}}\label{proof}
From now on, we assume that $k'$ is a splitting field for all groups considered below. The difference between Theorem \ref{main} and \cite[Theorem 3.1]{Marcus} is the coefficient rings. The ring $\O'$ in this note can be viewed as ``$\O$" in \cite{Marcus}. We will prove Theorem \ref{main} along Marcus' method.

Let $G:=S_n$ and $G^+:=A_n$. Let $b^+$ be a block of $\O G^+$ with a nontrivial abelian defect group $P$. Assume that $p=2$. Let $b'$ be a block of $\O'G^+$ such that $b'b^+\neq 0$. Then by \cite[Lemma 6.4]{Kessar_Linckelmann}, we can easily see that $P$ is also a defect group of $b'$. By \cite[Lemma (7.A)]{FH}, $P$ is a Klein four group. So in this case, Theorem \ref{main} follows from \cite[Theorem 3]{Huang}.

From now on, we assume that $p \geq3$. By Lemma \ref{block}, $b^+$ is also a block of $\O'G$ with defect group $P$.
From now on we can adopt the notation in Section \ref{notation}. The remaining proof consists of several steps.

By equation (\ref{2.1}), there are two situations. Assume that $b\neq {}^{\hat{\sigma}}b$, $b^+=b+{}^{\hat{\sigma}}b$. In this case, we have that $\O G^+b^+\cong\O Gb$ and $\O H^+c^+\cong\O Hc$.
By \cite[Theorem 7.6]{CR}, there is a splendid Rickard complex $X$ of $(\O Gb,\O Hc)$-bimodules. It follows that $X$ can be also viewed as a splendid Rickard complex of $(\O G^+b^+,\O H^+c^+)$-bimodules.

Now assume that $b^+=b={}^{\hat{\sigma}}b$, i.e., $b$ is {\it self-associated}. Then $c^+=c={}^{\hat{\sigma}}c$ and $\O Gb$ and $\O Hc$ are strongly $C_2$-graded algebras. We can apply \cite[Theorem 4.7]{Marcuscomm} if we showed that the splendid Rickard equivalence for $\O Gb$ and $\O Hc$ constructed in \cite[Theorem 7.6]{CR} is induced by a complex of $C_2$-graded bimodules. Since this equivalence is a composition of several splendid Morita or Rickard equivalences, we shall analyse each one of these equivalences.

Recall that the block $b$ corresponds uniquely to a $p$-core $\kappa$, and $P\cong C_p\times\cdots\times C_p$ ($w$ times) (see \cite[\S3]{CK}). Let $t:=n-pw$. Then by \cite[\S3]{CK}, the block algebra $\O Hc$ is isomorphic to $\O N_{S_{pw}}(P)\otimes_\O \O S_tc_0$, where $c_0$ is the block of defect zero of $\O S_t$ corresponding to the $p$-core $\kappa$. Recall that since $b$ is self-associated, $\kappa$ is also self-associated, namely that its diagram is symmetric with respect to the main diagonal. Since $N_{S_{pw}}(P)\cong N_{S_p}(C_p)\wr S_w\cong(C_p\rtimes C_{p-1})\wr S_w$, we have $\O N_{S_{pw}}(P)\cong \O((C_p\rtimes C_{p-1})\wr S_w)$.

Consider an abacus having $w+i(w-1)$ beads on the $i$-th runner, where $i=0,1,\cdots,p-1$, and let $\rho$ be the $p$-core have this abacus representation. Note that $\rho$ is self-associated. Let $r=|\rho|$, let $V$ be a set of cardinality $pw+r$, let $U_1, U_2,\cdots, U_w$ be disjoint subsets of $V$ of cardinality $p$, and let $U:=U_1\cup\cdots\cup U_w$. Let $e$ be a block of $\O S(V)$ with a defect group isomorphic to $P$ and corresponding to the $p$-core $\rho$. Here $S(V)$ is the symmetric group on $V$.
Let $\widetilde{N}$ be the subgroup of $S(U)$ consisting of permutations sending each $U_i$ into some $U_j$; we have that $\widetilde{N}$ is isomorphic to the wreath product $S_p\wr S_w$.
Let $N:=\widetilde{N}\times S(V\setminus U)$, and let $f\in \O N$ be the Brauer correspondent of $e$. Note that $N$ is a subgroup of $G$ containing $N_{S(V)}(P)$. By \cite[Lemma 5]{CK}, we have that the block algebra $\O Nf$ is isomorphic to $B_0(\O(S_p\wr S_w))\otimes_\O \O S_rf_0$, where $B_0(\O(S_p\wr S_w))$ is the principal block algebra of $\O(S_p\wr S_w)$, $f_0$ is the block of $\O S_r$ of defect zero corrsponding to the $p$-core $\rho$.

\medskip We consider the following pairs of block algebras:

\noindent(i) $\O Hc$ and $\O N_{S_{pw}}(P)\otimes_\O \O S_tc_0$;

\noindent(ii) $\O N_{S_{pw}}(P)\otimes_\O \O S_tc_0$ and $\O((C_p\rtimes C_{p-1})\wr S_w)\otimes_\O \O S_tc_0$;

\noindent(iii) $\O((C_p\rtimes C_{p-1})\wr S_w)\otimes_\O \O S_tc_0$ and $B_0(\O(S_p\wr S_w))\otimes_\O \O S_rf_0$;

\noindent(iv) $B_0(\O(S_p\wr S_w))\otimes_\O \O S_rf_0$ and $\O Nf$;

\noindent(v) $\O Nf$ and $\O S(V)e$;

\noindent(vi) $\O S(V)e$ and $\O Gb$.

If we showed that the algebras in each pair above are splendidly Morita or Rickard equivalent via a $C_2$-graded bimodule or complex of $C_2$-graded bimodules, then we showed that the splendid Rickard equivalence between $\O Gb$ and $\O Hc$ constructed in \cite[Theorem 7.6]{CR} are induced by a complex of $C_2$-graded bimodules.
By the discussion above, the algebras in the first pair are isomorphic; the algebras in the second pair are isomorphic; the algebras in the fourth pair are also isomorphic. So it suffices to consider the pairs (iii), (v) and (vi). If replacing the coefficient ring $\O$ with $\O'$, Marcus \cite{Marcus} have proved this, see \cite[3.4]{Marcus} for the pair (v), \cite[3.5 and 3.6]{Marcus} for the pair (iii) and \cite[3.7]{Marcus} for the pair (vi). We will prove the same thing for the ring $\O$ along Marcus' way.

\begin{lemma}
The block algebras $\O Nf$ and $\O S(V)e$ in the pair (v) above are splendidly Morita equivalent via a $C_2$-graded bimodule.
\end{lemma}

\noindent{\it Proof.} The $(\O'S(V)e, \O'S(V)e)$-bimodule $\O'S(V)e$ has vertex $\Delta P:=\{(u,u)|u\in P\}$. Let $M'$ be the Green correspondent of the $(\O'S(V),\O'S(V))$-bimodule $\O'S(V)e$ in $\O'(S(V)\times N)$. By \cite[page 179]{CK}, $M'$ induces a splendid Morita equivalence between $\O'S(V)e$ and $\O' Nf$.
By the definition of Green correspondent, $M'$ is the unique (up to isomorphism) direct summand of the $(\O'S(V)e,\O'Nf)$-bimodule $\O'S(V)e$ having $\Delta P$ as a vertex. Noting that $\O'S(V)e\cong \O'\otimes_\O \O S(V)e$, by \cite[Lemma 5.1]{Kessar_Linckelmann}, there is an indecomposable direct summand $M$ of the $(\O S(V)e,\O Nf)$-bimodule $\O S(V)e$, such that $M'\cong \O'\otimes_\O M$. By \cite[Proposition 4.5]{Kessar_Linckelmann}, $M$ induces a splendid Morita equivalence between $\O S(V)e$ and $\O Nf$.

The method to show that $M$ is a $C_2$-graded module is the same as Marcus' method, see the last paragraph of \cite[3.4]{Marcus}.~~$\hfill\square$

\medskip To prove that there is a $C_2$-graded Rickard equivalence between the block algebras in the third pais, by the first paragraph of \cite[3.5]{Marcus}, it suffices to show that there is a $C_2$-graded Rickard equivalence between $B_0(\O(S_p\wr S_w))$ and $\O((C_p\rtimes C_{p-1})\wr S_w)$, and a $C_2$-graded Rickard equivalence between $\O S_tc_0$ and $\O S_rf_0$.

\begin{lemma}\label{cyclic}
Let $b_0$ be the principal block of $\O S_p$, then there is a $C_2$-graded splendid Rickard equivalence between $\O S_pb_0$ and $\O (C_p\rtimes C_{p-1})$.
\end{lemma}

\noindent{\it Proof.} Since any $\O'$-algebra automorphism of $\O'S_p$ preserves the principle block, by (\ref{2.1}) and Lemma \ref{block}, $b_0$ is also the principal block idempotent of $\O'A_p$ and $\O A_p$. By \cite{Rou}, there is a splendid Rickard complex $X'$ for $\O'A_pb_0$ and $\O'(C_p\rtimes C_{\frac{{p - 1}}{2}})$. In \cite[3.5]{Marcus}, Marcus used \cite[Example 5.5]{Marcusonequivalences} to extend $X'$ to a $C_2$-graded splendid Rickard equivalence between $\O'S_pb_0$ and $\O'(C_p\rtimes C_{p-1})$. In \cite[Theorem 1.10]{Kessar_Linckelmann}, Kessar and Linckelmann showed that there is a splendid Rickard complex $X$ for $\O A_pb_0$ and $\O(C_p\rtimes C_{\frac{{p - 1}}{2}})$ such that $X'\cong \O'\otimes_\O X$. So we can prove the lemma by showing that $X$ extends to a $C_2$-graded splendid Rickard equivalence between $\O S_pb_0$ and $\O(C_p\rtimes C_{p-1})$.

The complexes $X'$ and $X$ have many similar properties, let us review the construction of them. Let $P$ be a Sylow $p$-subgroup of $A_p$. We identify $\O'(C_p\rtimes C_{\frac{{p - 1}}{2}})$ with $\O'N_{A_p}(P)$ and identify $\O(C_p\rtimes C_{p - 1})$ with $\O N_{S_p}(P)$. $X'$ is a two-term complex of the form
$$X':=\cdots\to 0\to X'_{-1}\xrightarrow{d'}X'_0\to 0\to \cdots,$$
where $X'_0$ is the $(\O'A_pb_0, \O'N_{A_p}(P))$-bimodule $\O'A_pb_0$. Decompose $X'_0$ as $M'\oplus U'$, where $M'$ is an indecomposable non-projective $(\O'A_pb_0, \O'N_{A_p}(P))$-bimodule and $U'$ is a projective $(\O'A_pb_0, \O'N_{A_p}(P))$-bimodule. Let $\pi':V'\to M'$ be a projective cover of $M'$. Then there is a certain direct summand $W'$ of $V'$, such that $X'_{-1}=W'\oplus U'$ and that $d'=\pi'_{W'}\oplus {\rm id}_{U'}$. To extend $X'$ to a $C_2$-graded splendid Rickard equivalence between $\O'S_pb_0$ and $\O'(C_p\rtimes C_{p-1})$, \cite[Example 5.5]{Marcusonequivalences} only used the following two properties: $M'$ is an indecomposable  non-projective summand of $X'_0$ (as an $(\O' A_pb_0, \O'N_{A_p}(P))$-bimodule), and $U'$ is a projective summand of $X'_0$; $W'$ and $V'/W'$ do not have any isomorphic direct summand.

By the proof of \cite[Theorem 1.10]{Kessar_Linckelmann}, $X$ is a two-term complex of the form
$$X:=\cdots\to 0\to X_{-1}\xrightarrow{d}X_0\to 0\to \cdots,$$
where $X_0$ is the $(\O A_pb_0, \O N_{A_p}(P))$-bimodule $\O A_pb_0$. Decompose $X_0$ as $M\oplus U$, where $M$ is an indecomposable non-projective $(\O A_pb_0, \O N_{A_p}(P))$-bimodule and $U$ is a projective $(\O A_pb_0, \O N_{A_p}(P))$-bimodule. Let $\pi:V\to M$ be a projective cover of $M$. Then there is a certain direct summand $W$ of $V$ satisfying $W'\cong \O'\otimes_{\O}W$, such that $X_{-1}=W\oplus U$ and that $d=\pi_{W}\oplus {\rm id}_{U}$. We still have the following two properties: $M$ is an indecomposable non-projective summand of $X_0$ (as an $(\O A_pb_0,\O N_{A_p}(P))$-bimodule), and $U$ is a projective summand of $X_0$; $W$ and $V/W$ do not have isomorphic direct summands.
We note that the blanket assumption in \cite{Marcusonequivalences} that $k'$ is ``big enough" was not used in most of arguments and conclusions in \cite{Marcusonequivalences}.
Hence the same argument in \cite[Example 5.5]{Marcusonequivalences} shows that $X$ extends to a $C_2$-graded splendid Rickard equivalence between $\O S_pb_0$ and $\O(C_p\rtimes C_{p-1})$. $\hfill\square$

\begin{lemma}
There is a $C_2$-graded Rickard equivalence between $B_0(\O(S_p\wr S_w))$ and $\O((C_p\rtimes C_{p-1})\wr S_w)$.
\end{lemma}

\noindent{\it Proof.} This follows from the same argument in the last paragraph of \cite[3.4]{Marcus}. $\hfill\square$

\medskip Next we prove that there is a $C_2$-graded Morita equivalence between the defect zero block algebras $\O S_tc_0$ and $\O S_rf_0$. It is clear that the integer $r$ above is larger than $1$, but the integer $t$ may equal to $0$ or $1$. Since the $p$-cores $\kappa$ and $\rho$ are self-associated, the idempotent $f_0$ decomposes as $f_0=f'+f''$ in $\O'A_r$; if $t\geq 2$, the idempotent $c_0$ decomposes as $c_0=c'+c''$ in $\O'A_r$ (see the second sentence of \cite[\S4]{FH} for the reason). Hence $c_0\in \O'A_t$, $f_0\in \O'A_r$. Since $k$ is a splitting field for symmetric groups, we also have that $c_0\in \O A_t$ and $f_0\in \O A_r$. If $t=0,1$, $\O'S_tc_0$ (resp. $\O S_tc_0$) is isomorphic to $\O'$ (resp. $\O$). Note that $t=0,1$ is possible, and this situation was not treated in Marcus' proof.

\begin{lemma}\label{dz}
There is a $C_2$-graded Morita equivalence between $\O S_tc_0$ and $\O S_r f_0$.
\end{lemma}

\noindent{\it Proof.} Since all blocks considered in this lemma is defect zero, so all modules of these algebras are projective. Since projective modules lifts uniquely (up to isomorphism) from $k$ to $\O$, we may replace $\O$ and $\O'$ by $k$ and $k'$, respectively. Since every finite group has a finite splitting field, we may assume that $k'$ is finite.

\medskip\noindent\textbf{3.4.1.} We first consider the more difficult case: $t\geq 2$. Keep the notation in the paragraph preceding this lemma. Let $T':=k'S_tc_0$, $T'_1:=k'A_tc_0$, $T:=k S_tc_0$, $R':=k'S_rf_0$, $R'_1:=k'A_rf_0$, and $R:=k S_rf_0$.  Let $x:=(1,2)$, an element in $S_t$ but not in $A_t$.
Let $y:=(1,2)$, an element in $S_r$ but not in $A_r$. We see that $c''=xc'x^{-1}$ and $f''=y^{-1}f'y$. Write $T'_x:=T'_1x$ and ${R'_y}^{\rm op}:=y{R'_1}^{\rm op}$, then $T'=T'_1\oplus T'_x$ and $R'^{\rm op}={R'_1}^{\rm op}\oplus {R'_y}^{\rm op}$ as $k'$-modules. With respect to these decompositions, $T'$ and $R'^{\rm op}$ are strongly $C_2$-graded algebras. Similarly, $T$ and $R^{\rm op}$ are strongly $C_2$-graded algebras.  Let $\Delta:=T'_1\otimes_k {R'_1}^{\rm op}\oplus T'_x\otimes_{k'}{R'_y}^{\rm op}$, then $\Delta$ is a $k'$-subalgebra of $T'\otimes_{k'} R'^{\rm op}$.

In \cite[3.6]{Marcus}, Marcus proved that there is a $C_2$-graded Morita equivalence between $T'$ and $R'$. Hence by Lemma \ref{gradcategory}, there is a left $(T'\otimes_{k'}R'^{\rm op})\rtimes \hat{C}_2$-module or a right $\hat{C}_2\ltimes(T'^{\rm op}\otimes_{k'} R')$-module $M'$ inducing a Morita equivalence between $T'$ and $R'$. We are going to show that there is a left $(T\otimes_kR^{\rm op})\rtimes \hat{C}_2$-module or a right $\hat{C}_2\ltimes(T^{\rm op}\otimes_{k} R)$-module $M$ such that $M'\cong k'\otimes_k M$. 

Let us review the construction of $M'$. Let $U$ be the unique (up to isomorphism) simple left $k'A_tc'$-module. $W$ be the unique (up to isomorphism) simple right $k'A_rf'$-module. The $x$-conjugate of $U$ is a $k'A_t$-module $\overline{U}$, which is isomorphic to $U$ as a $k'$-module. So we can write $\overline{U}$ as $\{\overline{u}~|~u\in U\}$, satisfying $\lambda\overline{u}=\overline{\lambda u}$ and $g\overline{u}:=\overline{(xgx^{-1})u}$ for any $\lambda\in k'$ and $g\in A_t$. We denote the $k'$-linear map $U\to \overline{U}, u\mapsto \overline{u}$ by $\ba_1$.
Similarly, the $y$-conjugate of $W$ can be written as $\overline{W}:=\{\overline{w}~|~w\in W\}$, satisfying $\lambda\overline{w}=\overline{\lambda w}$ and $\overline{w}h:=\overline{w(y^{-1}hy)}$ for any $\lambda\in k'$ and $h\in A_r$. We denote the $k'$-linear map $W\to \overline{W},w\mapsto \overline{w}$ by $\ba_2$. We see that $\overline{U}$ (resp. $\overline{W}$) is the unique (up to isomorphism) simple left $k'A_tc''$-module (resp. right $k'A_rf''$-module).

Let $V:=U\otimes_{k'} W\oplus \overline{U}\otimes_{k'} \overline{W}$, then $V$ is an $(T'_1\otimes_{k'}{R'_1}^{\rm op})$-module.
For any $u_1\in U$, $w_1\in W$, $\overline{u_2}\in \overline{U}$ and $\overline{w_2}\in \overline{W}$, define
$$(c_0x\otimes yf_0)(u_1\otimes w_1+\overline{u_2}\otimes\overline{w_2})=u_2\otimes w_2+\overline{u_1}\otimes\overline{w_1}.$$
One verifies that $V$ becomes a $\Delta$-module. By \cite[Theorem 3.4]{Marcusonequivalences}, the $(T',R')$-bimodule
$$M':=(T'\otimes_{k'}R'^{\rm op})\otimes_\Delta V$$
induces a $C_2$-graded Morita equivalence between $T'$ and $R'$. We can determine the structure of $M'$ more explicit. We define a $(T',R')$-bimodule $M''$: let
$$M'':= U\otimes_{k'} W\oplus \overline{U}\otimes_{k'} \overline{W}\oplus \overline{U}\otimes_{k'}W\oplus U\otimes_{k'}\overline{W}$$
as a $(T'_1,R'_1)$-bimodule. For any $u_1,u_2,u_3,u_4\in U$ and $w_1,w_2,w_3,w_4\in W$, define
\begin{equation}\label{a}
c_0x(u_1\otimes w_1+\overline{u_2}\otimes \overline{w_2}+\overline{u_3}\otimes w_3+u_4\otimes \overline{w_4})=u_3\otimes w_3+\overline{u_4}\otimes \overline{w_4}+\overline{u_1}\otimes w_1+u_2\otimes \overline{w_2}
\end{equation}
and define
\begin{equation}\label{b}
(u_1\otimes w_1+\overline{u_2}\otimes \overline{w_2}+\overline{u_3}\otimes w_3+u_4\otimes \overline{w_4})yf_0=u_4\otimes w_4+\overline{u_3}\otimes \overline{w_3}+\overline{u_2}\otimes w_2+u_1\otimes \overline{w_1}.
\end{equation}
Then one can verify that $M'\cong M''$ as $(T',R')$-bimodule. So we can identify $M'$ and $M''$.

If we decompose $T'\otimes_{k'} R'^{\rm op}$ as $T'_1\otimes_{k'} R'^{\rm op}\oplus T'_x\otimes_{k'} R'^{\rm op}$, then the $1$-component $M'_1$ and the $x$-component $M'_x$ of $M'$ can be chosen to be
\begin{equation}\label{componenta}
M'_1:=U\otimes_{k'} W \oplus \overline{U}\otimes_{k'}W~~{\rm and}~~M'_x:=\overline{U}\otimes_{k'} \overline{W}\oplus  U\otimes_{k'}\overline{W}
\end{equation}
(see \cite[Remarks 3.2 (b)]{Marcusonequivalences}).
If we decompose $T'^{\rm op}\otimes_{k'} R'$ as ${T'}^{\rm op}\otimes_{k'}R_1'\oplus {T'}^{\rm op}\otimes_{k'} R'_y$, then the $1$-component $M'_1$ and the $y$-component $M'_y$ of $M'$ can be chosen to be
\begin{equation}\label{componentb}
M'_1:=U\otimes_{k'} W \oplus U\otimes_{k'}\overline{W}~~{\rm and}~~M'_y:=\overline{U}\otimes_{k'} \overline{W}\oplus  \overline{U}\otimes_{k'}W
\end{equation}
By Lemma \ref{gradcategory}, $M'$ can be identified with a $(T'\otimes_{k'} R'^{\rm op})\rtimes \hat{C}_2$-module (resp. right $\hat{C}_2\ltimes(T'^{\rm op}\otimes_{k'} R')$-module) via
\begin{equation}\label{rho}
\hat{\rho}\cdot m_1=\hat{\rho}(1)m_1~~~{\rm and}~~~\hat{\rho}\cdot m_x=\hat{\rho}(x)m_x~~~({\rm resp.}~~~m_y\cdot\hat{\rho}=m_y\hat{\rho}(y))
\end{equation}
for any $\hat{\rho}\in \hat{C}_2$, $m_1\in M'_1$ and $m_x\in M'_x$ (resp. $m_y\in M'_y$).

Let $\Gamma:={\rm Gal}(k'/k)$, then $\Gamma$ is a cyclic group. Let $\tau$ be a generator of $\Gamma$. We claim  that ${}^\tau M'\cong M'$ as left $(T'\otimes_{k'} R'^{\rm op})\rtimes \hat{C}_2$-modules or as right $\hat{C}_2\ltimes(T'^{\rm op}\otimes_{k'} R')$-modules. It is easy to see that one of the four situations
\[\left\{ \begin{gathered}
  {}^\tau U\cong \overline{U} \hfill \\
  {}^\tau W\cong W \hfill \\
\end{gathered}  \right.;~~~\left\{ \begin{gathered}
  {}^\tau U\cong U \hfill \\
  {}^\tau W\cong \overline{W} \hfill \\
\end{gathered}  \right.;~~~\left\{ \begin{gathered}
  {}^\tau U\cong \overline{U}  \hfill \\
  {}^\tau W\cong \overline{W} \hfill \\
\end{gathered}  \right.;~~~\left\{ \begin{gathered}
  {}^\tau U\cong U \hfill \\
  {}^\tau W\cong W \hfill \\
\end{gathered}  \right.\]
should take place, where in each situation, the first isomorphism is an isomorphism of left $T'_1$-modules and  the second is an isomorphism of right $R'_1$-modules.
We are going to prove that in any of these situations, we can construct an isomorphism $\Pi: {}^\tau M'\cong M'$ of left $(T'\otimes_{k'} R'^{\rm op})\rtimes \hat{C}_2$-modules or right $\hat{C}_2\ltimes(T'^{\rm op}\otimes_{k'} R')$-modules.

Consider the first situation. Let $\varphi:{}^\tau U\cong \overline{U}$ be an isomorphism of left $T'_1$-modules and $\psi:{}^\tau W\cong W$ be an isomorphism of right $R'_1$-modules. Let
\begin{equation}\label{3.6}
\overline{\varphi}=\ba_1^{-1}\circ \varphi\circ\ba_1^{-1}: {}^\tau \overline{U}\to U
\end{equation}
and let
\begin{equation}\label{3.7}
\overline{\psi}=\ba_2\circ\psi\circ \ba_2^{-1}:{}^\tau \overline{W}\to \overline{W}.
\end{equation}
Using the property that $x$ (resp. $y$) is an element of order $2$, it is not difficult to check that $\overline{\varphi}$ (resp. $\overline{\psi}$) is an isomorphism of left $T'_1$-modules (resp. right $R'_1$-modules).  Define a map
$$\Pi: {}^\tau M'\to M'$$
$$u_1\otimes w_1+\overline{u_2}\otimes \overline{w_2}+\overline{u_3}\otimes w_3+u_4\otimes \overline{w_4}~\mapsto$$
$$\overline{\varphi}(\overline{u_3})\otimes \psi(w_3)+\varphi(u_4)\otimes\overline{\psi}(\overline{w_4})+\varphi(u_1)\otimes \psi(w_1)+\overline{\varphi}(\overline{u_2})\otimes \overline{\psi}(\overline{w_2}),$$
for any $u_1,u_2,u_3,u_4\in U$ and $w_1,w_2,w_3,w_4\in W$.
Clearly $\Pi$ is an isomorphism of $(T'_1,R'_1)$-bimodules. For any $m\in {}^\tau M'$ and $\hat{\rho}\in \hat{C}_2$,
the verifications of
$$\Pi(c_0xm)=c_0x\Pi(m),~~~\Pi(myf_0)=\Pi(m)yf_0~~~{\rm and}~~~\Pi(\hat{\rho}\cdot m)=\hat{\rho}\cdot\Pi(m)$$
are straightforward, by using equations (\ref{a}), (\ref{b}), (\ref{componenta}), (\ref{rho}), (\ref{3.6}) and (\ref{3.7}). Hence $\Pi$ is an isomorphism of $(T'\otimes_{k'} R'^{\rm op})\rtimes \hat{C}_2$-modules. The second situation is similar to the first situation.

Consider the third situation. Let $\varphi:{}^\tau U\cong \overline{U}$ be an isomorphism of left $T'_1$-modules and $\psi:{}^\tau W\cong \overline{W}$ be an isomorphism of right $R'_1$-modules. Let
\begin{equation}\label{3.8}
\overline{\varphi}:=\ba_1^{-1}\circ \varphi\circ\ba_1^{-1}: {}^\tau \overline{U}\to U
\end{equation}
and let
\begin{equation}\label{3.9}
\overline{\psi}:=\ba_2^{-1}\circ\psi\circ \ba_2^{-1}:{}^\tau \overline{W}\to W.
\end{equation}
Using the property that $x$ and $y$ are elements of order $2$, it is not difficult to check that $\overline{\varphi}$ (resp. $\overline{\psi}$) is an isomorphism of left $T'_1$-modules (resp. right $R'_1$-modules).  Define a map
$$\Pi: {}^\tau M'\to M'$$
$$u_1\otimes w_1+\overline{u_2}\otimes \overline{w_2}+\overline{u_3}\otimes w_3+u_4\otimes \overline{w_4}~\mapsto$$
$$\overline{\varphi}(\overline{u_2})\otimes \overline{\psi}(\overline{w_2})+\varphi(u_1)\otimes\psi(w_1)+\varphi(u_4)\otimes \overline{\psi}(\overline{w_4})+\overline{\varphi}(\overline{u_3})\otimes \psi(w_3),$$
for any $u_1,u_2,u_3,u_4\in U$ and $w_1,w_2,w_3,w_4\in W$.
Clearly $\Pi$ is an isomorphism of $(T'_1,R'_1)$-bimodules. For any $m\in {}^\tau M'$ and $\hat{\rho}\in \hat{C}_2$,
the verifications of
$$\Pi(c_0xm)=c_0x\Pi(m),~~~\Pi(myf_0)=\Pi(m)yf_0~~~{\rm and}~~~\Pi(m\cdot \hat{\rho})=\Pi(m)\cdot\hat{\rho}$$
are straightforward, by using equations (\ref{a}), (\ref{b}), (\ref{componentb}), (\ref{rho}), (\ref{3.8}) and (\ref{3.9}). Hence $\Pi$ is an isomorphism of right $\hat{C}_2\ltimes(T'^{\rm op}\otimes_{k'} R')$-modules. The last situation is similar to the third situation.

Since all blocks in the lemma are defect zero and since $p>2$, it is easy to see that $M'$ is a projective  $(T'\otimes_{k'} R'^{\rm op})\rtimes \hat{C}_2$-module or right $\hat{C}_2\ltimes(T'^{\rm op}\otimes_{k'} R')$-module. By \cite[Lemma 6.2 (c)]{Kessar_Linckelmann}, that there is a $(T\otimes_kR^{\rm op})\rtimes \hat{C}_2$-module or right $\hat{C}_2\ltimes(T^{\rm op}\otimes_{k} R)$-module $M$ such that $M'\cong k'\otimes_k M$. Then by \cite[Proposition 4.5 (c)]{Kessar_Linckelmann}, $M$ induces a $C_2$-graded Morita equivalence between $T$ and $R$. Now we have completed the proof for $t\geq 2$.

\medskip\noindent\textbf{3.4.2.} Let us turn to the easier case: $t=0,1$. Now $T'=k'S_tc_0=k'$, $T=kS_0f_0=k$. We regard $T'$ as a $C_2$-graded $k'$-algebra by setting $T'_1:=k'A_tf_0=k'$ and $T'_y:=0$. Here we identify $C_2$ to the subgroup $\{1,y\}$ of $S_r$ for notational convenience. Similarly, $T$ is $C_2$-graded $k$-algebra where $T_1:=kA_tf_0=k$ and $T_y:=0$.

Since $W$ is a simple right $k'A_rf'$-module, it is easy to see that $W\otimes_{k'A_r} k'S_r$ is a simple right $k'S_rf_0$-module. We can also regard $W\otimes_{k'A_r} k'S_r$ as a left $T'$-module, where the element in $S_t$ acts trivially on this module. The $(T',R')$-bimodule $M':=W\otimes_{k'A_r} k'S_r$ induces a splendid Morita equivalence between $T'$ and $R'$. The module structure of $M'$ can be described more explicit:
$M'=W\oplus \overline{W}$ as a $(T',R'_1)$-bimodule (where $S_t$ acts trivially); for any $w_1,w_2\in W$,
$(w_1+\overline{w_2})yf_0=w_2+\overline{w_1}$.
It is also easy to see that the Morita equivalent is also a $C_2$-graded Morita equivalence between $T'$ and $R'$, where the $1$-component $M_1':=W$ and the $y$-component $M'_y:=\overline{W}$.

As in the first case, we complete the proof by showing that ${}^\tau M'\cong M'$ as $\hat{C}_2\ltimes (T'^{\rm op}\otimes_k R')$-modules. If there is an isomorphism $\psi: {}^\tau W\cong W$ of right $R_1'$-modules, take $\overline{\psi}$ to be the map defined as in (\ref{3.7}). If there is an isomorphism $\psi: {}^\tau W\cong \overline{W}$ of right $R_1'$-modules, then take $\overline{\psi}$ to be the map defined as in (\ref{3.9}).
Define a map
$$\Pi:{}^\tau M'\to M',~~w_1+\overline{w_2}\mapsto \varphi(w_1)+\overline{\psi}(\overline{w_2}),$$
for any $w_1,w_2\in W$.
Then it is straightforward to check that $\Pi$ is an isomorphism of $\hat{C}_2\ltimes (T'^{\rm op}\otimes_k R')$-modules.   $\hfill\square$

\medskip It remains to show that there is a $C_2$-graded Rickard equivalence between the block algebras $\O Gb$ and $\O S(V)e$ in the pair (vi). In \cite{CR}, Chuang and Rouquier proved that there is a splendid Rickard equivalence between them. Marcus proved that over the large enough coefficient ring $\O'$, the Rickard
equivalence can be lifted to a $C_2$-graded Rickard equivalence (see \cite[3.7]{Marcus}). \cite[Lemma 10.2.13]{Rou} plays a key role in Marcus' proof. But \cite[Lemma 10.2.13]{Rou} requires the coefficient ring to have algebraically closed residue field. The following Lemma is a variation of \cite[Lemma 10.2.13]{Rou}, which will be used in our proof for the ring $\O$.

\begin{lemma}\label{extension}
Let $G_1$ be a normal subgroup of $G_2$ with $E=G_2/G_1$ a cyclic $p'$-group. Let $M$ be a $G_2$-stable $\O G_1$-module. Let $g\in G_2/G_1$ generating $E$ and let $\varphi\in {\rm End}_\O (M)$ such that
$\varphi(g^{-1}hg(m))=h\varphi(m)$ for all $h\in G_2$ and $m\in M$. Assume that there exists $\alpha\in {\rm End}_{\O G_1}(M)$ such that $\alpha^{|E|}=\varphi^{|E|}g^{-|E|}$ as elements in ${\rm End}_\O(M)$. Then there exists an $\O G_2$ module $\widetilde{M}$ extending $M$.
\end{lemma}

\noindent{\it Proof.} This follows by the same arguments as in the proof of \cite[Lemma 10.2.13]{Rou}. In \cite[Lemma 10.2.13]{Rou}, the assumption that the residue field is algebraically closed ensures the existence of $\alpha$. $\hfill\square$

\begin{lemma}
There is a $C_2$-graded Morita $C_2$-graded Rickard equivalence between $\O Gb$ and $\O S(V)e$ in the pair (vi).
\end{lemma}

\noindent{\it Proof.} According to Marcus's description, the Rickard equivalence between $\O S_nb$ and $\O S(V)e$ is obtained as a composition of equivalences between blocks forming a so-called $[w:m]$ pair, defined as follows. Assume that $a$ is a block of weight $w$ of $\O S_n$ corresponding to an abacus whose $j$-th runner has $m$ more beads than the $(j-1)$-th runner. Switching the number of beads on these two runners, one obtains a block $a'$ of weight $w$ of $\O S_{n-m}$.

\medskip\noindent\textbf{3.6.1.} If $m\geq w$, Scopes \cite{Scopes} proved that $\O S_na$ and $\O S_{n-m}a'$ are Morita equivalent. We treat this situation first. Observe that $M:=a\O S_n a'$ is an $(\O S_na, \O S_{n-m}a'\otimes_\O \O S_m)$-bimodule. Here, we view $S_{n-m}$ (resp. $S_m$) as a subgroup of $S_n$ via the first $n-m$ letters (resp. the last $m$ letters). Then the Morita equivalence is induced by the $(\O S_na, \O S_{n-m}a')$-bimodule $M\otimes_{\O S_m}\O$. One of the problems is that usually $a$ and $a'$ are not both self-associated, we need to discuss four situations:

Case 1. Both $a$ and $a'$ are self-associated, namely that $a={}^{\hat{\sigma}}a=a^+$, $a'={}^{\hat{\sigma}}a'=a'^+$ (see (\ref{2.1})). Denote by $x=(1,2)$, an element in $S_n$ but not in $A_n$. We identify $C_2$ to $\{1,x\}$. In this case, both $A:=\O S_n a$ and $B:=\O S_{n-m}a'$ are $C_2$-graded algebras, where $A_1=\O A_na$, $A_{x}=\O A_nxa$, $B_1=\O A_{n-m}a'$, and $B_x=\O A_{n-m}xa'$. We can decompose $M:=a\O S_n a'$ as $M=M_1\oplus M_x$ where $M_1=a\O A_n a'$, $M_x=a\O A_nx a'$. With this decomposition, $M$ is a $C_2$-graded $(A,B)$-bimodule.

Case 2. Both $a$ and $a'$ are not self-associated, namely that $a\neq{}^{\hat{\sigma}}a$, $a'\neq{}^{\hat{\sigma}}a'$, $a^+=a+{}^{\hat{\sigma}}a$, $a'^+=a'+{}^{\hat{\sigma}}a'$ (see (\ref{2.1})). In this case, $\O S_na\cong \O A_na^+$ and $\O S_{n-m}a'\cong \O A_{n-m}a'^+$. Let $A:=\O A_n a^+$ and $B:=\O A_{n-m}a'^+$. Both $A$, $B$ are $C_2$-graded algebras, where $A_1=\O A_na^+$, $A_{x}=0$, $B_1=\O A_{n-m}a'^+$, and $B_x=0$. We can decompose $M:=a\O S_n a'$ as $M=M_1\oplus M_x$ where $M_1=a\O S_n a'$, $M_x=0$. With this decomposition, $M$ is a $C_2$-graded $(A,B)$-bimodule.

Case 3. $a$ is self-associated but $a'$ is not. Let $A:=\O S_n a$ and $B:=\O A_{n-m}a'^+$. Let $y:=(1,n)$, an element in $S_n$, commutating with all elements in $S_{n-m}$, but not in $A_n$.  Both $A$, $B$ are $C_2$-graded algebras, where $A_1=\O A_na$, $A_{x}=\O A_nya$, $B_1=\O A_{n-m}a'^+$, and $B_x=0$.  We can decompose $M:=a\O S_n a'$ as $M=M_1\oplus M_x$ where $M_1=a\O A_n a'$, $M_x=a\O A_nya'$. With this decomposition, $M$ is a $C_2$-graded $(A,B)$-bimodule.

Case 4. $a'$ is self-associated but $a$ is not. Let $A:=\O A_n a^+$ and $B:=\O S_{n-m}a$. Both $A$, $B$ are $C_2$-graded algebras, where $A_1=\O A_na^+$, $A_{x}=0$, $B_1=\O A_{n-m}a'$, and $B_x=\O A_{n-m}xa'$.
Consider the $\hat{\sigma}$-conjugate of the $(A,B)$-bimodule of $M:=a\O S_n a'$. One verifies that the map
$$\varphi:~M\to {}^{\hat{\sigma}}M,~~u\mapsto {}^{\hat{\sigma}}u$$
is an isomorphism of $(A,B)$-bimodules (noting that the elements in $M$ can be viewed as elements in $A$, see (\ref{action}) for the definition of ${}^{\hat{\sigma}}u$). So $M$ is $\hat{C}_2$-stable. Clearly that $\varphi^2={\rm Id}_M$ as an element in ${\rm End}_\O (M)$. If we want to use Lemma \ref{extension} to show that $M$ extends to $(A\otimes_\O B^{\rm op})\rtimes \hat{C}_2$, we need to show that there exists $\alpha\in {\rm End}_{A\otimes_\O B^{\rm op}}(M)$ such that $\alpha^2={\rm Id}_M$. Define a map $\alpha: M\to M$ sending $u$ to $uy$, then $\alpha$ is a desired map.

Actually, the argument in Case 4 also applies to the first three cases, but the arguments in the first three cases make the structure of the $C_2$-graded bimodule more explicit.  By the discussion above, in any case, $M$ is an $(A\otimes_\O B^{\rm op})\rtimes \hat{C}_2$-module. Moreover, $M$ is an $((A\otimes_\O B^{\rm op})\rtimes \hat{C}_2,\O S_m)$-bimodule. Consequently, $M\otimes_{\O S_m}\O$ is an $(A\otimes_\O B^{\rm op})\rtimes \hat{C}_2$-module, hence a $C_2$-graded $(A,B)$-bimodule.

\medskip\noindent\textbf{3.6.2.} For arbitrary $m$, the Rickard complex is a generalization of Scopes' bimodule. In this situation, the method to prove that the Rickard complex extends to a complex of $C_2$-graded bimodules is the same as Marcus' method in \cite[3.7.2]{Marcus}. We only need to replace \cite[3.7.1]{Marcus} there with 3.6.1.  $\hfill\square$

\medskip Now we have proved that the block algebras in each pair listed at the beginning of this section are splendidly Morita or Rickard equivalent via a $C_2$-graded bimodule or complex of $C_2$-graded bimodules. So we have shown that the splendid Rickard equivalence between $\O Gb$ and $\O Hc$ constructed in \cite[Theorem 7.6]{CR} are induced by a complex of $C_2$-graded bimodules. By \cite[Theorem 4.7]{Marcuscomm}, Theorem \ref{main} holds.

\bigskip\noindent\textbf{Acknowledgements.}\quad I wish to express gratitude to Prof. Joseph Chuang for suggesting me to consider the descent question on alternating groups and guiding me to some references; to Prof. Radha Kessar and Prof. Markus Linckelmann for many helpful discussions and for their kindly hosting me during my one-year visit at City, University of London; to Prof. Andrei Marcus for answering my questions on his papers.
I also thank China Scholarship Council for financial support (No. 202006770016).


\end{document}